\documentclass[leqno,11pt]{amsart}
\usepackage{amssymb,amsmath}

\calclayout

\newtheorem{theorem}{Theorem}

\newtheorem{lemma}[theorem]{Lemma}

\newcommand{\cF}{{\mathcal F}}
\newcommand{\cH}{{\mathcal H}}

\newcommand{\mC}{{\mathbb C}}
\newcommand{\mF}{{\mathbb F}}
\newcommand{\mP}{{\mathbb P}}
\newcommand{\mQ}{{\mathbb Q}}
\newcommand{\mZ}{{\mathbb Z}}

\newcommand{\fg}{{\mathfrak g}}
\newcommand{\fh}{{\mathfrak h}}
\newcommand{\fs}{{\mathfrak s}}
\newcommand{\ft}{{\mathfrak t}}

\newcommand{\cl}{\operatorname{cl}}
\newcommand{\diag}{\operatorname{diag}}
\newcommand{\Fr}{\operatorname{Fr}}
\newcommand{\gr}{\operatorname{gr}}
\newcommand{\Int}{\operatorname{Int}}
\newcommand{\red}{\operatorname{red}}
\newcommand{\Sym}{\operatorname{Sym}}

\title{The virtual Poincar\'e polynomials of homogeneous spaces}
\author{Michel~Brion}
\address{Universit\'e de Grenoble I\\
D\'epartement de Math\'ematiques\\
Institut Fourier, UMR 5582 du CNRS\\
38402 Saint-Martin d'H\`eres Cedex, France}
\email{Michel.Brion@ujf-grenoble.fr, Emmanuel.Peyre@ujf-grenoble.fr}

\author{Emmanuel~Peyre}
\address{}
\email{}

\date{}
\begin{document}

\begin{abstract}
We factor the virtual Poincar\'e polynomial of every homogeneous space 
$G/H$, where $G$ is a complex connected linear algebraic group and $H$
is an algebraic subgroup, as $t^{2u} (t^2-1)^r Q_{G/H}(t^2)$ 
for a polynomial $Q_{G/H}$ with non-negative integer coefficients. 
Moreover, we show that $Q_{G/H}(t^2)$ divides the virtual Poincar\'e
polynomial of every regular embedding of $G/H$, if $H$ is
connected.
\end{abstract}

\maketitle

\section*{Introduction and statement of the results}

One associates to every complex algebraic variety $X$  (possibly
singular, or reducible) its {\sl virtual Poincar\'e polynomial}
$P_X(t)$, uniquely determined by the following properties:

\smallskip

\noindent
(i) (additivity) $P_X(t) = P_Y(t) + P_{X-Y}(t)$ for every closed
subvariety $Y$.

\smallskip

\noindent
(ii) If $X$ is smooth and complete, then 
$P_X(t)=\sum_m \dim H^m(X) \; t^m$ 
is the usual Poincar\'e polynomial.

\smallskip

\noindent
Then $P_X(t) = P_Y(t) \; P_F(t)$ for every fibration $F\to X\to Y$
which is locally trivial for the Zariski topology. 

\smallskip

Specifically, we have
$$
P_X(t)=\sum_{j,m} (-1)^{j+m} \dim\gr_W^m(H^j_c(X))\; t^m,
$$
where $gr_W^m(H^j_c(X))$ denotes the $m$-th subquotient of the weight
filtration on the $j$-th cohomology group of $X$ with compact supports
and complex coefficients (see \cite{F1} 4.5 and \cite{DL}). More
generally, the mixed Hodge structure on $H^*_c(X)$ yields a polynomial
$E_X(s,t)$ in two variables, satisfying the same properties of
additivity and multiplicativity, and such that $P_X(t)=E_X(-t,-t)$ 
(see \cite{DK} and \cite{BD} \S 3 for more details).

\smallskip

In this paper, we investigate the $E$-polynomials of
homogeneous spaces under linear algebraic groups, and of their 
{\sl regular} embeddings in the sense of \cite{BDP}. It turns out that
these polynomials behave much better than the usual Poincar\'e
polynomials; the latter are generally unknown for homogeneous
spaces. To state our main results, we introduce the following
notation.

\smallskip

Let $G$ be a complex connected linear algebraic group and let $H$ be a
closed subgroup. Let $r_H$ (resp.~$u_H$) be the rank (resp.~the
dimension of a maximal unipotent subgroup) of $H$, and define
similarly $r_G$, $u_G$. Choose maximal reductive subgroups 
$H^{\red}\subseteq H$, $G^{\red}\subseteq G$ such that
$H^{\red}\subseteq G^{\red}$, and maximal tori $T_H\subseteq H^{\red}$,
$T_G=T\subseteq G^{\red}$ such that $T_H\subseteq T$; let $W_H$, $W_G=W$
be the corresponding Weyl groups. The Lie algebras of $G$, $H,\ldots$
will be denoted $\fg$, $\fh,\ldots$ 

\smallskip

The group $W_H$ acts on the Lie algebra $\ft_H$ and on its ring of
polynomial functions, $\mC[\ft_H]=R(T_H)$. The invariant subring
$\mC[\ft_H]^{W_H}=R(H)$ is a finitely generated, graded algebra over
$\mC$, isomorphic to $\mC[\fh^{\red}]^{H^{\red}}$. Its Hilbert series
$\sum_{m=0}^{\infty} \dim R(H)_m\;t^m$ is the expansion
of a rational function of $t$, denoted $F_H(t)$.

\smallskip

Since $G$ is connected, $R(G)$ is a polynomial ring, and there
exists a graded subspace $\cH$ of $R(T)$ such that the multiplication
map induces an isomorphism of $R(G)\otimes\cH$ onto $R(T)$. Moreover,
$\cH$ is isomorphic to the cohomology space of the flag variety
$\cF(G)$, with complex coefficients. This isomorphism doubles degrees,
and the Hodge structure on $H^*(\cF(G))$ is pure. Therefore,
the Poincar\'e polynomial $P_{\cF(G)}$ is even, and we have
$$
E_{\cF(G)}(s,t)=P_{\cF(G)}((st)^{1/2}) \text{ and }
\frac{1}{(1-t)^{r_G}}=F_T(t)=F_G(t) P_{\cF(G)}(t^{1/2}).
$$
Moreover, we have
$$
P_{\cF(G)}(q^{1/2}) = \vert\cF(G)(\mF_q)\vert
$$ 
for every finite field $\mF_q$ with $q$ elements. Here
$\vert\cF(G)(\mF_q)\vert$ denotes the number of points over $\mF_q$ of 
$\cF(G)$ regarded as the flag variety of the split $\mZ$-form of
$G^{\red}$. 

\smallskip

Our first main result generalizes this to an arbitrary homogeneous
space $G/H$, with some twists. Notice that both $G$ and its closed
subgroup $H$ are defined over a finitely generated subring of $\mC$,
so that $(G/H)(\mF_q)$ makes sense for a large power $q$ of a large
prime number.

\begin{theorem}\label{1}
(a) With preceding notation, the virtual Poincar\'e polynomial $P_{G/H}$
is even, and we have
$$
E_{G/H}(s,t)=P_{G/H}((st)^{1/2}) \text{ and }
F_H(t)=F_G(t)\;t^{\dim(G/H)} P_{G/H}(t^{-1/2}).
$$
Moreover, we have for all large $q$:
$$
\vert(G/H)(\mF_q)\vert=P_{G/H}(q^{1/2}).
$$

\noindent
(b) There exists a polynomial $Q_{G/H}$ with non-negative
integer coefficients, such that
$$
P_{G/H}(t^{1/2})=t^{u_G-u_H} (t-1)^{r_G-r_H} Q_{G/H}(t).
$$
Moreover,
$$
Q_{G/H}(t)=Q_{G^{\red}/H^{\red}}(t).
$$
The degree of $Q_{G/H}$ equals $\dim\cF(G)-\dim\cF(H^0)$, with
leading coefficient $1$, and $Q_{G/H}(1)$ equals 
$\frac{\vert W_G\vert}{\vert W_H\vert}$. 

\noindent
(c) If $H$ is connected, then 
$$
Q_{G/H}(t)=\frac{P_{\cF(G)}(t^{1/2})}{P_{\cF(H)}(t^{1/2})}
=t^{\dim \cF(G)-\dim\cF(H)}Q_{G/H}(t^{-1}).
$$
In particular, $Q_{G/H}(0)=1$.
\end{theorem}

It follows that $u_G-u_H$, $r_G-r_H$ and $Q_{G/H}$ depend only on the
complex algebraic variety $G/H$ (in fact, $r_G-r_H$ is a topological
invariant, see \cite{AP} 4.3).

\smallskip

As another consequence, the Poincar\'e polynomial of the flag variety
of a semi-simple group is divisible by the Poincar\'e polynomial of
the flag variety of every semi-simple subgroup, and the quotient has
non-negative coefficients. 

\smallskip

Theorem \ref{1} is proved in Section 1 by arguments of equivariant
cohomology ; it would be interesting to deduce it from a deeper
motivic result. Notice that (a) can be deduced from the fibration 
$$
G/H\to BH\to BG,
$$ 
where $BH$ (resp.~$BG$) denotes the classifying space of $H$
(resp.~$G$); then the cohomology ring of $BH$ is isomorphic to $R(H)$
with degrees doubled, so that the Poincar\'e series of $BH$ is
$F_H(t^2)$. If moreover $H$ is connected, then 
$$
P_{G/H}(t^{1/2})=\frac{P_G(t^{1/2})}{P_H(t^{1/2})}
=t^{u_G-u_H} (t-1)^{r_G-r_H}
\frac{P_{\cF(G)}(t^{1/2})}{P_{\cF(H)}(t^{1/2})},
$$
as follows from \cite{DL} Theorem 6.1 (ii); and a similar relation
holds for $\vert(G/H)(\mF_q)\vert$, by Lang's theorem. 

So the main point of Theorem \ref{1} is (b), especially the
non-negativity of coefficients of $Q_{G/H}$. We deduce it (together
with (a) and (c)) from a geometric construction that may be of
independent interest. In loose words, we obtain a locally trivial
fibration (for the Zariski topology)
$$
S\to G/H\to Z
$$ 
where $S$ is a torus of dimension $r_G-r_H$, and $Z$ is an algebraic
variety satisfying Poincar\'e duality and whose cohomology is purely
algebraic (see Lemmas \ref{torus} and \ref{pure} for a precise
statement). Thus, $E_{G/H}(s,t)=(1-st)^{r_G-r_H} E_Z(s,t)$, and
$E_Z(s,t)$ is the value at $(st)^{1/2}$ of the Poincar\'e polynomial
of $H^*_c(Z)$. In the case where $G$ and $H$ have the same rank, it
follows that $P_{G/H}(t)$ is the Poincar\'e polynomial of
$H^*_c(G/H)$.

\medskip

Next we turn to the $E$-polynomials of {\sl regular embeddings}. 
Recall from \cite{BDP} that a regular embedding of $G/H$ is a smooth
complex algebraic variety $X$ endowed with an algebraic action of $G$,
such that:

\noindent
(i) $X$ contains an open orbit isomorphic to $G/H$.

\noindent
(ii) The complement of this open orbit is a union of smooth
irreducible divisors (the {\sl boundary divisors}), with normal
crossings.

\noindent
(iii) Every orbit closure is a partial intersection of the
boundary divisors, and its normal bundle contains an open orbit.

\smallskip

Recall also that those homogeneous spaces under a connected reductive 
group $G$ which admit a complete regular embedding are exactly the
{\it spherical} homogeneous spaces, i.e., those where a Borel subgroup
of $G$ acts with an open orbit. 

\smallskip

Since every regular embedding $X$ contains only finitely many orbits,
we have 
$$
E_X(s,t)=P_X((st)^{1/2})
$$ 
by Theorem \ref{1} and additivity. Therefore, it suffices to consider
the virtual Poincar\'e polynomial $P_X$. Our second main result yields
a factorization of that polynomial:

\begin{theorem}\label{2}
Let $X$ be a regular embedding of $G/H$, where $H$ is connected. Then,
for every orbit $G/H'$ in $X$, the polynomial $Q_{G/H}(t)$ divides
$Q_{G/H'}(t)$, and the quotient has non-negative integer coefficients.

As a consequence, there exists a polynomial $R_X(t)$ with integer
coefficients, such that
$$
P_X(t^{1/2}) = Q_{G/H}(t) R_X(t).
$$

If moreover $X$ is complete, then the coefficients of $R_X(t)$ are
non-negative. 
\end{theorem}

The assumption that $H$ is connected cannot be suppressed, as shown by
an example at the end of Section 2. This section is devoted to the
proof of Theorem 2. Again, the main point is the non-negativity of
coefficients of $R_X(t)$; for this, we show that the equivariant
cohomology ring of $X$ is a free module of finite rank over a
polynomial subring generated by $R(H)$ and indeterminates of degree
$2$. It would be interesting to obtain a topological interpretation of
the polynomial $R_X(t)$. However, the factorization 
$P_X(t^{1/2}) = Q_{G/H}(t) R_X(t)$ does not originate in a fibration
with total space $X$, as shown by the following simple example.

\smallskip

Consider the complex projective space $X=\mP^{2m+1}$ of odd dimension,
where the projective special orthogonal group 
$G=SO(2m+2)/\{\pm 1\}$ acts linearly. Then $X$ consists of $2$ orbits:
the quadric $Q^{2m}$, and its complement with isotropy group 
$H\cong O(2m+1)/\{\pm 1\}\cong SO(2m+1)$, a connected subgroup;
one checks that $X$ is a regular completion of $G/H$. We have
$$
P_{G/H}(t^{1/2})=P_{\mP^{2m+1}}(t^{1/2})-P_{Q^{2m}}(t^{1/2})
=t^m(t^{m+1}-1),
$$
so that $Q_{G/H}(t)=t^m+t^{m-1}+\cdots+1$ and that $R_X(t)=t^{m+1}+1$. 
How to explain the factorization
$$
P_{\mP^{2m+1}}(t^{1/2})=t^{2m+1}+t^{2m}+\cdots+1
=(t^m+t^{m-1}+\cdots+1)(t^{m+1}+1)
$$
in topological terms ?

\smallskip

Notice that the complex projective space $\mP^{2m}$ of even dimension is
a regular completion of the homogeneous space $SO(2m+1)/O(2m)$ (where
$O(2m)$ is not connected) by the quadric $\mQ^{2m-1}$; this yields
$Q_{SO(2m+1)/O(2m)}(t)=1$. 

\smallskip

These are examples of complete symmetric varieties. In fact, the
Poincar\'e polynomials of all such varieties were determined by 
De Concini and Springer (see \cite{DS}) who deduced the virtual
Poincar\'e polynomials of adjoint symmetric spaces. Their results 
were the starting point for the present work, as the factorizations of
Theorems 1 and 2 can be seen on examples of \cite{DS}. 

\smallskip

For instance, by Theorem 2, the virtual Poincar\'e polynomial of any
regular embeding $X$ of a connected reductive group $G$ (viewed as a
homogeneous space under the action of $G\times G$ by left and right
multiplication) is divisible by $Q_G(t^2)=P_{\cF(G)}(t^2)$. When $G$
is semi-simple adjoint and $X$ is its canonical completion, this
agrees with the closed formula for $P_X(t)$ given in \cite{DS} p. 96.

\section{Proof of Theorem 1}

In what follows, we use \cite{EZ} as a general reference for mixed
Hodge structure, and \cite{S} for algebraic groups. 

We begin with an easy reduction to the case where both groups $G$ and
$H$ are reductive. Let $R_u(H)$ be the unipotent radical of $H$. This
unipotent group is isomorphic, as an algebraic variety, to some
$\mC^u$. Since $H$ is the semi-direct product of $R_u(H)$ with
$H^{\red}$, we have $u=u_H-u_{H^{\red}}$. The quotient map $G\to G/H$
factors through
$$
p:G/H^{\red}\to G/H,
$$
a fibration with fiber $R_u(H)\cong \mC^u$. Thus, the pullback map
$H^*(G/H^{\red})\to H^*(G/H)$ is an isomorphism of mixed Hodge
structures. By Poincar\'e duality, it follows that
$$
E_{G/H^{\red}}(s,t) = (st)^u \; E_{G/H}(s,t).
$$

We now show that
$$
\vert(G/H^{\red})(\mF_q)\vert = q^u \; \vert(G/H)(\mF_q)\vert
$$
for $q$ such that $H^{\red}$ is defined over $\mF_q$ and that $H$ is
the semidirect product of $R_u(H)$ with $H^{\red}$ over
$\overline{\mF_q}$. This follows from Grothendieck's trace formula; 
as an alternative proof using elementary arguments of Galois
descent, we check that
$$
\pi:(G/H^{\red})(\mF_q)\to (G/H)(\mF_q)
$$
is surjective with all fibers of order $q^u$. We denote 
$\Fr_q$ the Frobenius endomorphism of $G(\overline{\mF_q})$, with
fixed point subgroup $G(\mF_q)$.

Let $x\in G(\overline{\mF_q})$ such that $xH\in (G/H)(\mF_q)$. Then
$x^{-1}\Fr_q(x)\in H(\overline{\mF_q})$. Thus, we can write 
$x^{-1}\Fr_q(x)=yz$ where $y\in R_u(H)(\overline{\mF_q})$ and
$z\in H^{\red}(\overline{\mF_q})$. Since  
$R_u(H)(\overline{\mF_q})$ is connected and invariant under
$\Int(z)\circ \Fr_q$, there exists
$h\in R_u(H)(\overline{\mF_q})$ such 
that $y=h z \Fr_q(h^{-1}) z^{-1}$. Thus, 
$x^{-1}\Fr_q(x)=h z \Fr_q(h^{-1})$. 
Replacing $x$ by $xh$, we may assume that 
$x^{-1}\Fr_q(x)\in H^{\red}(\overline{\mF_q})$. This proves the
surjectivity of $\pi$.

Let now $x,y \in G(\overline{\mF_q})$ such that 
$xH^{\red},yH^{\red}\in (G/H^{\red})(\mF_q)$ and that $y\in xH$. We
may assume that $y=xz$ where $z\in R_u(H)(\overline{\mF_q})$. Then
$H^{\red}(\overline{\mF_q})$ contains 
$x^{-1}\Fr_q(x)$ and 
$y^{-1}\Fr_q(y)=z^{-1}x^{-1}\Fr_q(x)\Fr_q(z)$. Since 
$H^{\red}(\overline{\mF_q})$ normalizes $R_u(H)(\overline{\mF_q})$ and
their intersection is trivial, it follows that
$z^{-1}x^{-1}\Fr_q(x)\Fr_q(z)\Fr_q(x^{-1})x=1$.
Therefore, $xzx^{-1}\in (xR_u(H)x^{-1})(\mF_q)$, and
$xR_u(H)x^{-1}$ is a $\Fr_q$-stable connected unipotent
group of dimension $u$. So every fiber of $\pi$ has order $q^u$.

Therefore, if Theorem \ref{1} holds for $G/H^{\red}$, then it holds
for $G/H$, and 
$$
Q_{G/H^{\red}}(t)=Q_{G/H}(t).
$$
So we may assume that $H=H^{\red}$. Then, using the fibration
$$
G/H\to G/R_u(G)H\cong G^{\red}/H^{\red}
$$
with fiber $R_u(G)$, one reduces similarly to the case where
$G=G^{\red}$.

\medskip

We assume from now on that $G$ and $H$ are reductive; as a consequence,
$G/H$ is affine.

\begin{lemma}\label{torus}
The following conditions are equivalent for a subtorus $S$ of $T$, 
with Lie algebra $\fs \subseteq \ft$:

\noindent
(i) All isotropy subgroups of $S$ acting on $G/H$ are
finite, and $S$ is maximal for this property.

\noindent
(ii) $\fs\oplus w\ft_H=\ft$ for all $w\in W$.

As a consequence, there exist subtori $S$ satisfying (i), and all of
them have dimension $r_G-r_H$. Moreover, the double coset space
$S\backslash G/H$ is an affine algebraic variety, with at worst
quotient singularities by finite abelian groups.
\end{lemma}

\begin{proof}
Let $g\in G$, then the finiteness of the isotropy group of $gH$ in $S$
is equivalent to: $\fs\cap{\rm Ad}(g)\fh=0$. As there are only
finitely many isotropy groups for a torus action on an algebraic
variety, the finiteness of all isotropy groups for the $S$-action on
$G/H$ is equivalent to: $\fs\cap{\rm Ad}(G)\fh=0$. Since
$$
\fs\cap{\rm Ad}(G)\fh = \fs\cap(\ft\cap{\rm Ad}(G)\ft_H)
= \fs\cap W\ft_H,
$$
this amounts to: $\fs\cap w\ft_H=\{0\}$ for all $w\in W$.

Now $\ft$ has a $W$-invariant rational structure, defined by the
lattice of differentials at $1$ of one-parameter subgroups of $T$; the
rational subspaces are exactly the Lie algebras of subtori. Moreover,
any rational subspace $\fs$ intersecting trivially all subspaces
$w\ft_H$ is contained in a rational complement to all these
subspaces. This proves equivalence of conditions (i) and (ii), and the
assertion on existence of subtori $S$ and their dimension. For any
such subtorus $S$, all orbits in the affine variety $G/H$ are closed,
and the isotropy groups are finite abelian groups. This implies the
latter assertion.
\end{proof}

\medskip

\noindent
{\sl Remark.} Lemma \ref{torus} extends to arbitrary homogeneous
spaces $G/H$, except for the assertion that $S\backslash G/H$ is an
affine algebraic variety. In fact, the quotient space 
$S\backslash G/H$ may well be non-separated if $G/H$ is not
affine. For example, let $G={\rm SL}(2)$ and let $H$ be its standard
unipotent subgroup. The diagonal torus $D\cong\mC^*$ of $G$ acts on
$G/H\cong\mC^2 - \{0\}$  by $t\cdot (x,y)=(tx,t^{-1}y)$. All isotropy
groups are trivial, but the quotient space is a classical example of
a non-separated scheme : the affine line with its origin doubled.

\medskip

Next choose a subtorus $S$ of $T$ satisfying the conditions of
Lemma \ref{torus} and let 
$$
Z=S\backslash G/H
$$ 
with quotient map $f:G/H\to Z$. Then there exists a decomposition of 
$Z$ into finitely many disjoint, locally closed
subvarieties $Z_j$ ($j\in J$), together with finite subgroups $F_j$
($j\in J$) of $S$, such that every $f^{-1}(Z_j)$ is equivariantly
isomorphic to $S/F_j \times Z_j$. Since $S/F_j$ is a torus of
dimension $r_G-r_H$, we have $E_{S/F_j}(s,t)=(st-1)^{r_G-r_H}$, whence
$$
E_{G/H}(s,t)=(st-1)^{r_G-r_H} E_{Z}(s,t).
$$
Likewise, we have for all large $q$:
$$
\vert(G/H)(\mF_q)\vert=(q-1)^{r_G-r_H}\vert Z(\mF_q)\vert.
$$

Since $Z$ has at worst finite quotient singularities, it satisfies
Poincar\'e duality over $\mC$. As a consequence, each closed algebraic
subvariety of codimension (say) $r$ in $Z$ has a cohomology class in
$H^{2r}(Z)$. This yields the (degree doubling) cycle map
$$
\cl:A^*(Z)\to H^*(Z),
$$ 
where the left hand side is the Chow group of $Z$,
graded by codimension (see \cite{F2} Chapter 19).

\begin{lemma}\label{pure} With preceding notation, $\cl$ is an
isomorphism over $\mC$. Moreover, the graded ring 
$H^*(Z)$ is isomorphic to $R(S)\otimes_{R(G)} R(H)$,
and the usual Poincar\'e polynomial of $Z$ equals
$$
\frac{F_S(t^2)F_H(t^2)}{F_G(t^2)}=
\frac{F_H(t^2)}{(1-t^2)^{r_G-r_H}F_G(t^2)}.
$$
\end{lemma}

\begin{proof}
We use equivariant cohomology, see e.g. \cite{H}. Consider the 
action of $T$ on $G/H$, then the equivariant cohomology ring
$H^*_T(G/H)$ is clearly isomorphic to $H^*_H(G/T)$. Since
$H^*_G(G/T)=H^*(BT)=R(T)$ is a free module of rank $\vert W\vert$ over 
$H^*_G(pt)=H^*(BG)=R(G)$, the Eilenberg-Moore spectral sequence (see
\cite{H} III.2) yields an isomorphism
$$
H^*_H(G/T)\cong H^*(BH)\otimes_{H^*(BG)} H^*_G(G/T),
$$
that is,
$$
H^*_T(G/H)\cong R(T)\otimes_{R(G)} R(H).
$$
This is a commutative, positively graded algebra, finite and free of
rank $\vert W\vert$ over its subring $R(H)$. The latter is a
Cohen-Macaulay ring of dimension $r_H$. Thus, the ring $H^*_T(G/H)$ is
Cohen-Macaulay of dimension $r_H$ as well, with Poincar\'e series
$$
\frac{F_T(t^2)F_H(t^2)}{F_G(t^2)}=
\frac{F_H(t^2)}{(1-t^2)^{r_G}F_G(t^2)}.
$$
Since the subtorus $S$ of $T$ acts on $G/H$ with finite isotropy
groups, we have
$$
H^*_T(G/H) \cong H^*_{T/S}(S\backslash G/H) \cong H^*_{T/S}(Z).
$$
This is a finitely generated module over $H^*_{T/S}(pt)=R(T/S)$. But
$T/S$ is a torus of dimension $r_H$, so that $R(T/S)$ is a polynomial
ring in $r_H$ variables. Since $H^*_T(G/H)$ is Cohen-Macaulay of
dimension $r_H$ and finite over $R(T/S)$, it is a free module over
that ring, by the Auslander-Buchsbaum formula (see \cite{E} 19.3). 
By the Eilenberg-Moore spectral sequence again, 
it follows that the canonical map
$$
\mC\otimes_{R(T/S)}H^*_{T/S}(Z)\to H^*(Z)
$$
is an isomorphism. Therefore, we have
$$
H^*(Z) \cong 
\mC\otimes_{R(T/S)}R(T)\otimes_{R(G)} R(H).
$$
But $\mC\otimes_{R(T/S)}R(T)\cong R(S)$; thus, we obtain 
$H^*(Z)\cong R(S)\otimes_{R(G)} R(H)$. 
Moreover, $H^*(Z)$ is the quotient of $H^*_{T/S}(Z)$ by
a regular sequence consisting of $r_H$ homogeneous elements of degree
$2$. Therefore, the usual Poincar\'e polynomial of $Z$
equals
$$
(1-t^2)^{r_H}\frac{F_T(t^2)F_H(t^2)}{F_G(t^2)}=
\frac{F_H(t^2)}{(1-t^2)^{r_G-r_H}F_G(t^2)}.
$$

It remains to compare cohomology of $Z$ with its Chow group. For this,
we use equivariant intersection theory, see \cite{EG} and also
\cite{B}. The equivariant Chow group with complex coefficients (graded
by codimension) $A^*_T(G/H)_{\mC}$ is again isomorphic to
$R(T)\otimes_{R(G)} R(H)$, by \cite{B} Corollary 12. Moreover, for any
scheme $X$ with an action of $T$, the natural map
$$
R(S)\otimes_{R(T)} A^*_T(X)_{\mC}\to A^*_S(X)_{\mC}
$$
is an isomorphism (to see this, one reduces to the case where the
quotient $X\to X/T$ exists and is a principal $T$-bundle, and one
argues as in \cite{B}, p. 17). As a consequence, the map
$$
R(S)\otimes_{R(G)} R(H)\to A^*_S(G/H)_{\mC}
$$
is an isomorphism; it follows that the cycle map
$$
\cl :A^*_S(G/H)_{\mC}\to H^*_S(G/H)=H^*(Z)
$$
defined in \cite{EG} 2.8, is an isomorphism as well. Finally,
$A^*_S(G/H)_{\mC}\cong A^*(S\backslash G/H)_{\mC}=A^*(Z)_{\mC}$ 
by \cite{EG} Proposition 4
and Theorem 4.
\end{proof}

\medskip

\noindent
{\sl Remark}. By Lemma \ref{pure}, the Betti numbers of 
$Z=S\backslash G/H$ are independent of the choice of $S$. But the
algebra structure of $H^*(Z)$ may depend on $S$, as shown by the
example where $H={\rm SL}(2)\times{\rm SL}(2)$ is embedded diagonally
in $H\times H=G$. Furthermore, there may exist no subtorus $S$ acting
on $G/H$ with finite constant isotropy groups; this happens, for
instance, if $G={\rm SL}(3)$ and $H={\rm SO}(3)$.

\medskip

As a final preparation for the proof of Theorem 1, we need the
following easy result of invariant theory.

\begin{lemma}\label{series}
We have 
$$
\lim_{t\to 1} (1-t)^{r_H} F_H(t) = \frac{1}{\vert W_H\vert}.
$$
Moroever, the degree of the rational function $F_H(t)$ is at most
$-\dim\cF(H^0)$, with equality if $H$ is connected.
\end{lemma}

\begin{proof}
The former assertion is a (well-known) consequence of Molien's formula
for the invariant ring $R(H)=\mC[\ft_H]^{W_H}$:
$$
F_H(t)=\frac{1}{\vert W_H\vert}\sum_{w\in W_H}
\frac{1}{\det_{\ft_H}(1-tw^{-1})}.
$$

For the latter assertion, recall that $R(H^0)$ is a graded polynomial
ring with homogeneous generators of degrees $d_1\leq\cdots\leq d_r$,
where $r=r_H$. Thus, the degree of $F_{H^0}(t)$ is
$-d_1-\cdots-d_r=-\dim\cF(H^0)$. Moreover, denoting $\Gamma$ the
finite group $H/H^0$, we have an exact sequence
$$
1\to W_{H^0}\to W_H\to \Gamma\to 1.
$$
Thus, $\Gamma$ acts on $R(H^0)$ with invariant subring $R(H)$. Since
$R(H^0)$ is a graded polynomial ring, it contains a graded
$\Gamma$-stable subspace $V$ such that the map $\Sym(V)\to R(H^0)$
is an isomorphism. It follows that $V$ decomposes as a direct sum of
homogeneous components $V_d$ ; the increasing sequence of their
degrees (with multiplicities given by the dimensions of the $V_d$) is
the same as $(d_1,\ldots,d_r)$. Now
$$
F_H(t)=\frac{1}{\vert\Gamma\vert}\sum_{\gamma\in\Gamma}
\frac{1}{\prod_d \det_{V_d}(1-t^d\gamma^{-1})}
$$ 
is a sum of rational functions of the same degree, equal to
$-d_1-\cdots-d_r=-\dim\cF(H^0)$. 
\end{proof}

We can now complete the proof of Theorem 1. By Lemma \ref{pure}, the
cohomology of $Z$ vanishes in all odd degrees, and every space
$H^{2m}(Z)$ is generated by algebraic classes. Thus, the Hodge
structure on that space is pure of type $(m,m)$, and the same holds
for the dual space $H^{2m}_c(Z)$. In other words,
$$
E_Z(s,t)=\sum_m \dim H^{2m}_c(Z)\; (st)^m.
$$
Using Poincar\'e duality and Lemma \ref{pure}, it follows that
$$
E_Z(s,t)=(st)^{\dim(Z)}
\frac{F_H((st)^{-1})}{(1-(st)^{-1})^{r_G-r_H}F_G((st)^{-1})},
$$
so that
$$
E_{G/H}(s,t)=\frac{(st)^{\dim(G/H)} F_H((st)^{-1})} {F_G((st)^{-1})}.
$$
On the other hand, we have
$$
\vert Z(\mF_q)\vert = \sum_m \dim H^{2m}_c(Z)\; q^m.
$$
For, by Grothendieck's trace formula \cite{Se}, one has
$$
\vert Z(\mF_q)\vert = \sum_m(-1)^m\text{Tr}(\Fr_q,
H^m_c(Z_{\overline \mF_q},\mQ_l))
$$
the equality, for large $q$, then follows from the proper base change
theorem and the fact that the cycle class map is an isomorphism.  
(alternatively, one may show directly that 
$$
\vert (G/H)(\mF_q)\vert = 
\frac{q^{\dim(G/H)} F_H(q^{-1})}{F_G(q^{-1})},
$$
by arguments of Galois descent). This implies (a). 
Taking degrees in the equality of rational functions
$$
F_H(t)=F_G(t)t^{\dim(G/H)} \; P_{G/H}(t^{-1/2})
$$
and using Lemma \ref{series}, we obtain that $P_{G/H}(t^{1/2})$ is
divisible by 
$$
t^{\dim(G/H)-\dim \cF(G)+\dim\cF(H)}=t^{u_G-u_H}.
$$
Thus, we can write 
$P_{G/H}(t^{1/2})=t^{u_G-u_H} (t-1)^{r_G-r_H}Q_{G/H}(t)$ 
for a polynomial $Q_{G/H}(t)$ with integer coefficients. Since
$P_Z(t)=t^{u_G-u_H} Q_{G/H}(t)$, these coefficients are non-negative. 
Moreover, Lemma \ref{series} implies that 
$Q_{G/H}(1)=\frac{\vert W_G\vert}{\vert W_H\vert}$.

For any irreducible variety $X$, the degree of $P_X(t)$ is $2\dim(X)$,
with leading coefficient $1$. It follows that the degree of
$Q_{G/H}(t)$ is $\dim\cF(G)-\dim\cF(H^0)$, with leading coefficient
$1$. This completes the proof of (b). Finally, (c) follows from (a),
(b) and Poincar\'e duality for $\cF(G)$ and $\cF(H)$.

\section{Proof of Theorem 2}

Let $Y$ be an orbit in $X$. Replacing $X$ by the union of all
orbits whose closure contains $Y$ (an open $G$-invariant subset of 
$X$), we may assume that $Y$ is closed in $X$. Then $Y$ is the
transversal intersection of boundary divisors, say $X_1,\ldots,X_r$. 
Choose $x\in Y$ and denote by $H'$ its isotropy subgroup. Then $H'$
acts on the normal space to $Y$ at $x$; this action is diagonalizable
and given by $r$ linearly independent characters, see \cite{BDP}. This
defines a surjective group homomorphism $H'\to(\mC^*)^r$, whence an
exact sequence
$$
1\to K\to H'\to (\mC^*)^r\to 1
$$
where $K$ is the kernel of the $H'$-action on the normal space.
Let $K^{\red}$ be a maximal reductive subgroup of $K$.

We claim that $K^{\red}$ is contained in a conjugate of $H$. To check
this, consider the linear action on $K^{\red}$ on the tangent space
$T_x X$ and choose a $K^{\red}$-invariant complement $N$ to the
$K^{\red}$-invariant subspace $T_x Y$; by construction, $K^{\red}$
fixes $N$ pointwise. Then we can choose a $K^{\red}$-invariant
subvariety $Z$ of $X$, such that $Z$ is smooth at $x$ and that 
$T_x Z=N$. Therefore, $K^{\red}$ fixes pointwise a neighborhood of $x$
in $Z$, and this neighborhood meets the open orbit $G/H$.

Thus, we may assume that $K^{\red}$ is contained in $H$. Since $H$ is
connected, we can apply \cite{DL} Theorem 6.1 (ii) to the fibration
$G/K^{\red}\to G/H$ with fiber $H/K^{\red}$, to obtain
$$
P_{G/K^{\red}}(t)=P_{G/H}(t) P_{H/K^{\red}}(t).
$$
Together with Theorem \ref{1}, it follows that
$$
Q_{G/K}(t)=Q_{G/H}(t) Q_{H/K^{\red}}(t).
$$
On the other hand, the right action of
$H'/K\cong (\mC^*)^r$ on $G/K$ defines a principal $(\mC^*)^r$-bundle
$G/K\to G/H'$. All such bundles are locally trivial, whence
$P_{G/K}(t)=(t^2-1)^r P_{G/H'}(t)$, and 
$$
Q_{G/K}(t)=Q_{G/H'}(t).
$$ 
So, $Q_{G/H}(t)$ divides $Q_{G/H'}(t)$ and the quotient has
non-negative coefficients. 

By additivity, it follows that $Q_{G/H}(t)$ divides
$P_X(t^{1/2})$ ; the quotient is an even polynomial, $R_X(t)$. Since
$Q_{G/H}(0)=1$, the coefficients of $R_X(t)$ are integers. However,
their non-negativity for complete $X$ is not an obvious fact, because of
the factor $t^{u_G-u_{H'}}(t-1)^{r_G-r_{H'}}$ in each
$P_{G/H'}(t^{1/2})$. For this reason, we shall present an alternative
proof of the existence of $R_X(t)$, which will also yield this
non-negativity property. 

We begin by relating the virtual Poincar\'e polynomial
$P_X(t)$ to equivariant cohomology of $X$. If $V$ is a $\mZ$-graded
complex vector space such that every homogeneous component $V_m$ is
finite dimensional, let
$F_V(t)=\sum_{m=-\infty}^{\infty} \dim(V_m)\; t^m$ 
be its Poincar\'e series. If $X$ is a variety where $G$ acts
algebraically, then $H^*_G(X)$ is a finitely generated, graded module
over $H^*(BG)=R(G)$. As a consequence, the series $F_{H^*_G(X)}(t)$
is the expansion of a rational function, for which we use the same
notation.

\begin{lemma}\label{equiv}
For every regular embedding $X$, the rational function
$F_{H^*_G(X)}(t)$ is even, and 
$$
F_{H^*_G(X)}(t^{1/2}) = F_G(t)\; t^{\dim(X)} P_X(t^{-1/2}).
$$
\end{lemma}

\begin{proof}
In the case where $X=G/H$ is a unique orbit, we have 
$H^*_G(X)\cong H^*(BH)\cong R(H)$, whence
$F_{H^*_G(X)}(t)=F_H(t^2)$. So the assertion follows from
Theorem 1.

In the general case, choose a closed orbit $Y$ in $X$, of codimension
$r$, with complement $U$. The inclusion map $i:Y\to X$ defines a Gysin
morphism
$$
i_*:H^*_G(Y)\to H^*_G(X),
$$ 
of degree $2r$. By \cite{BDP}, this map and the restriction map
$H^*_G(X)\to H^*_G(U)$ fit into a short exact sequence 
$$
0\to H^*_G(Y)\to H^*_G(X)\to H^*_G(U)\to 0.
$$
It follows that
$$
F_{H^*_G(X)}(t)=t^{2r}F_{H^*_G(Y)}(t)+F_{H^*_G(U)}(t).
$$
Since $P_X=P_Y+P_U$, our assertion follows by induction.
\end{proof}

\medskip

\noindent
{\sl Remark}. Lemma \ref{equiv} admits a simpler formulation in terms
of equivariant Borel-Moore homology $H_*^G(X)$, as defined in
\cite{EG}. Indeed, by Poincar\'e duality, the rational function
$F_{H_*^G(X)}(t)$ is even, and
$$
F_{H_*^G(X)}(t^{1/2}) = F_G(t^{-1}) P_X(t^{1/2}).
$$
In fact this holds, more generally, for every variety $X$ where $G$
acts with finitely many orbits.
\medskip

Next let $X_1,\ldots,X_n$ be the boundary divisors of the regular
embedding $X$, and let $z_1,\ldots,z_n\in H^2_G(X)$ be their
equivariant cohomology classes. In the ring $H^*_G(X)$, consider the
ideal $I_X$ of $H^*_G(X)$ generated by $z_1,\ldots,z_n$, and the ideal
$J_X$, kernel of the restriction map
$$
\rho:H^*_G(X)\to H^*_G(G/H)\cong R(H).
$$
Clearly, $I_X$ is contained in $J_X$, and the latter ideal is
prime. Moreover, $\rho$ is surjective by \cite{BDP}, so that we have
an exact sequence
$$
0\to J_X\to H^*_G(X)\to R(H)\to 0.
$$
Examples show that $I_X$ may differ from $J_X$; but these ideals are
closely related, as shown by the following result.

\begin{lemma}\label{power}
We have $J_X^{2^N}\subseteq I_X$, where $N$ denotes the number of
$G$-orbits in $X$.
\end{lemma}

\begin{proof}
We argue by induction on $N$. If $N=1$, then $X=G/H$ so that both
$I_X$ and $J_X$ are trivial. In the general case, we use the notation
of the proof of Lemma \ref{equiv}. The (surjective) restriction map 
$H^*_G(X)\to H^*_G(U)$ sends $I_X$ (resp.~$J_X$) onto $I_U$
(resp.~$J_U$).

Let $\alpha\in J_X$. Since $J_U^{2^{N-1}}\subseteq I_U$ by the
induction assumption, we may assume that 
$$
\alpha^{2^{N-1}}=i_*\beta
$$ 
for some $\beta\in H^*_G(Y)$. Now we have in $H^*_G(X)$:
$$
\alpha^{2^N}=(i_*\beta)\cup(i_*\beta)=i_*(\beta\cup i^*i_*\beta)
=i_*(\beta^2\cup i^*i_*1)=(i_*\beta^2)\cup(i_*1),
$$
by the projection formula. Moreover, $i_*1$ is the equivariant
cohomology class of $Y$ in $X$. Since $Y$ is a transversal
intersection of $r$ boundary divisors, say $X_1,\ldots,X_r$, we have
$i_*1=z_1\cdots z_r\in I_X$, and $\alpha^{2^N}\in I_X$ as well.
\end{proof}

Since $H$ is connected, $R(H)$ is a graded polynomial ring, so that we
can choose a graded subalgebra $R$ of $H^*_G(X)$ that restricts
isomorphically to $H^*_G(G/H)\cong R(H)$ via $\rho$.  

\begin{lemma}\label{finite}
$H^*_G(X)$ is finite over its subring generated by $R$ and
$z_1,\ldots,z_n$.
\end{lemma}

\begin{proof}
Since the algebra $H^*_G(X)$ is positively graded, it suffices to
prove that the quotient 
$$
H^*_G(X)/(z_1,\ldots,z_n)=H^*_G(X)/I_X
$$ 
is a finitely generated $R$-module. By Lemma \ref{power},
$H^*_G(X)/I_X$ is a quotient of $H^*_G(X)/J_X^{m}$ for some positive
integer $m$. Consider the finite filtration of $H^*_G(X)/J_X^m$ by the
powers of the image of $J_X$, and notice that all the subquotients
$J_X^pH^*_G(X)/J_X^{p+1}H^*_G(X)$
are finite modules over $H^*_G(X)/J_X=R(H)$. Since the latter is
isomorphic to $R$, the assertion follows.
\end{proof}

We now need the following variant of the Noether normalization theorem.

\begin{lemma}\label{noether}
Let $A$ be a finitely generated, positively graded algebra over an
infinite field $k$. Let $y_1,\ldots,y_m$ be homogeneous, algebraically
independent elements of $A$ and let $z_1,\ldots,z_n$ be homogeneous
elements of degree $1$, such that $A$ is finite over its
subalgebra generated by 
$y_1,\ldots,y_m,z_1,\ldots,z_n$. Then there
exist a non-negative integer $n'$ and homogeneous elements
$y'_1,\ldots,y'_m,z'_1,\ldots,z'_{n'}$ of $A$ such that:

\noindent
(i) $y'_i-y_i \in k[z_1,\ldots,z_n]$ for $1\leq i\leq m$.

\noindent
(ii) $z'_1,\ldots,z'_{n'}$ are linear combinations of
$z_1,\ldots,z_n$. 

\noindent
(iii) $y'_1,\ldots,y'_m,z'_1,\ldots,z'_{n'}$ are algebraically
independent, and $A$ is finite over the subring that they generate.
\end{lemma}

\begin{proof} The argument is similar to that of the classical Noether
normalization theorem, see \cite{E} 13.1; we present it for
completeness. We argue by induction on $n$, the case where $n=0$ being
trivial. In the general case, we may assume that
$y_1,\ldots,y_m,z_1,\ldots,z_n$ are algebraically dependent, and we
choose a polynomial relation
$$
P(y_1,\ldots,y_m,z_1,\ldots,z_n)=0.
$$
We may assume that this relation is homogeneous and
involves $z_n$. Let $d_1,\ldots,d_m$ be the degrees of
$y_1,\ldots,y_m$. Define $y'_1,\ldots,y'_m,z'_1,\ldots,z'_{n-1}$ by
$$
y_i=y'_i + a_i z_n^{d_i} ,~z_j=z'_j+b_jz_n
$$
where $a_1,\ldots,a_m,b_1,\ldots,b_{n-1}$ are in $k$. Then
$$
P(y'_1 + a_1 z_n^{d_1}, \ldots, y'_m + a_m z_n^{d_m},
z'_1 + b_1 z_n, \ldots, z'_{n-1}+b_{n-1}z_n,z_n)=0.
$$
Regarding the right-hand side as a polynomial in $z_n$, the
coefficient of the leading term equals
$P(a_1,\ldots,a_m,b_1,\ldots,b_{n-1},1)$. Since $k$ is infinite and by
our assumptions on $P$, we may choose
$a_1,\ldots,a_m,b_1,\ldots,b_{n-1}$ so that this coefficient is 
non-zero. Then $z_n$ is integral over the subring $A'$ of $A$
generated by $y'_1,\ldots,y'_m$ and $z'_1,\ldots,z'_{n-1}$. We conclude
by the induction assumption for $A'$.
\end{proof}

We can now show that $Q_{G/H}(t)$ divides $P_X(t^{1/2})$. Apply Lemma
\ref{noether} to the algebra $H^*_G(X)$ and to homogeneous,
algebraically independent generators of its polynomial subalgebra $R$;
then we obtain another polynomial subalgebra $R'$ (restricting
isomorphically to $R(H)$) and linear combinations
$z'_1,\ldots,z'_{n'}$ of $z_1,\ldots,,z_n$, such that $H^*_G(X)$ is
finite over its polynomial subring
$R'[z'_1,\ldots,z'_{n'}]$. Let $f(t)$ be the associated Hilbert
polynomial, then
$$
F_{H^*_G(X)}(t^{1/2})=\frac{F_H(t) f(t)}{(1-t)^{n'}}.
$$
Moreeover, $f(1)$ is the rank of the $R'[z'_1,\ldots,z'_{n'}]$-module
$H^*_G(X)$, a positive integer. On the other hand, we 
have by Lemma \ref{equiv}:
$$
F_{H^*_G(X)}(t^{1/2}) = F_G(t)\;t^{\dim(G/H)} P_X(t^{-1/2})
$$
and, by Theorem \ref{1}:
$$
F_H(t) = F_G(t) \; t^{\dim(G/H)-u_G+u_H}
(t^{-1}-1)^{r_G-r_H} Q_{G/H}(t^{-1}).
$$
This yields
$$
P_X(t^{1/2}) = t^{n'+u_G-u_H} (t-1)^{r_G-r_H-n'} Q_{G/H}(t) f(t^{-1}).
$$
Since $f(1)Q_{G/H}(1)\neq 0$, we must have $r_G-r_H-n'\geq 0$; and
since $Q_{G/H}(0)=1$, the Laurent polynomial 
$t^{n'+u_G-u_H}(t-1)^{r_G-r_H-n'}f(t^{-1})$ must be a
polynomial. Thus, $Q_{G/H}(t)$ divides $P_X(t^{1/2})$. 

If moreover $X$ is complete, then the $R(G)$-module $H^*_G(X)$ is
free by \cite{BDP}. Thus, the ring $H^*_G(X)$ is Cohen-Macaulay of
dimension $r_G$. Since this ring is finite over 
$R'[z'_1,\ldots,z'_{n'}]$, a polynomial subring, $H^*_G(X)$ is a
free module over that subring, and we have $r_G=r_H+n'$. Therefore,
the Hilbert polynomial $f(t)$ has non-negative coefficients, so that
the same holds for the polynomial
$$
t^{n'+u_G-u_H}f(t^{-1})=\frac{P_X(t^{1/2})}{Q_{G/H}(t)}.
$$

\medskip

\noindent
{\sl Example.} We show that Theorem 2 does not extend to all
homogeneous spaces $G/H$. Let $G={\rm SL}(2)\times{\rm SL}(2)$ with
maximal torus $T=D\times D$, where $D$ denotes the diagonal torus of
${\rm SL}(2)$. Let  
$n=\left(\begin{matrix}0&1\cr -1&0\cr\end{matrix}\right)$, then the
element $(n,n)$ of $G$ normalizes $T$.
Let $H$ be the subgroup of $G$ generated by $T$ and by $(n,n)$. 
The homogeneous space $G/H$ is spherical, and we have
$T_H=T$. Denoting by $x$, $y$ the obvious coordinates on $\ft$,
one obtains $R(G)=\mC[x^2,y^2]$ and $R(H)=\mC[x^2,xy,y^2]$, whence
$$
F_G(t)=\frac{1}{(1-t^2)^2},~F_H(t)=\frac{1+t^2}{(1-t^2)^2},~ 
P_{G/H}(t^{1/2})=t^4+t^2 \text{ and }Q_{G/H}(t)=1+t^2.
$$

We now construct a regular completion $X$ of $G/H$, such that
$P_X(t^{1/2})$ is not divisible by $Q_{G/H}(t)$. Consider the variety 
$Y=\mP^1\times\mP^1\times\mP^1\times\mP^1$ 
where $G$ acts by
$(g_1,g_2)(a,b,c,d)=(g_1 a, g_1 b, g_2 c, g_2 d)$.
Then $Y$ is a regular embedding of $G/T$. Moreover, the right action
of $(n,n)$ on $G/T$ extends to the involution $\sigma$ of $Y$,
defined by
$\sigma(a,b,c,d)=(b,a,d,c)$.
The fixed point subset $Y^{\sigma}$ is the closed
$G$-orbit, $\diag(\mP^1)\times\diag(\mP^1)$.
Since the actions of $G$ and $\sigma$ commute, $G$ acts on the
quotient $Y/\sigma$. The latter is singular along the image $Z$ of
$Y^{\sigma}$; the normal space to $Y/\sigma$ at every point of $Z$ is 
isomorphic to the quotient of $\mC^2$ by the involution 
$(s,t)\mapsto (-s,-t)$. Thus, blowing up $Z$ along $Y/\sigma$ yields 
a smooth projective embedding $X$ of $G/H$.

One may check that $X$ is regular and that
$P_X(t^{1/2})= t^4 + 3t^3 + 6t^2 + 3t + 1$,
which is prime to $Q_{G/H}(t)=t^2+1$. One may also check that
$Q_{G/H'}(t)$ equals $t+1$ or $(t+1)^2$ for the other orbits; thus,
$Q_{G/H}(t)$ is prime to all other $Q_{G/H'}(t)$.


\begin{thebibliography}{100}

\bibitem{AP} C.~Allday and V.~Puppe: Cohomological methods in
transformation groups, Cambridge Studies in Advanced Mathematics 
{\bf 32} (1993), Cambridge University Press, 1993.

\bibitem{BD} V.~V.~Batyrev and D.~I.~Dais: Strong McKay
correspondence, string-theoretic Hodge numbers and mirror symmetry,
{\sl Topology} {\bf 35} (1996), 901-929. 

\bibitem{BDP} E.~Bifet, C.~De Concini and C.~Procesi: Cohomology of
regular embeddings, {\sl Adv. Math.} {\bf 82} (1990), 1-34.

\bibitem{B} M.~Brion: Equivariant cohomology and equivariant
intersection theory, pp. 1-37 in: Representation theories and
algebraic geometry, Nato ASI Series {\bf 514}, Kluwer, 1998.

\bibitem{DK} V.~I.~Danilov and A.~A.~Khovanskii: Newton polyhedra
and an algorithm for computing Hodge-Deligne numbers, 
{\sl Math. USSR Izvestiya} {\bf 29} (1987), 279-298.

\bibitem{DS} C.~De Concini and T.~A.~Springer: Betti numbers of
complete symmetric varieties, pp. 87-107 in: Geometry Today, Progress
in Math. {\bf 60}, Birkh\"auser, 1985.

\bibitem{DL} A.~Dimca and G.~Lehrer: Purity and equivariant weight
polynomials, pp. 161-181 in: Algebraic groups and Lie groups, Australian
Math. Soc. Lecture Series {\bf 9}, Cambridge University Press, 1997.

\bibitem{EG} D.~Edidin and W.~Graham: Equivariant intersection theory,
{\sl Invent. math.} {\bf 131} (1998), 595-634. 

\bibitem{E} D.~Eisenbud: Commutative algebra with a view towards
algebraic geometry, Graduate Text in Math. {\bf 150}, Springer-Verlag,
1995.

\bibitem{EZ} F.~El Zein: Introduction \`a la th\'eorie de Hodge mixte, 
Hermann, 1991.

\bibitem{F1} W.~Fulton: Introduction to Toric Varieties, Annals of
Math. Studies {\bf 131}, Princeton University Press, 1993.

\bibitem{F2} W.~Fulton: Intersection Theory, Ergeb. Math.
{\bf 2}, Springer-Verlag 1998.

\bibitem{H} W.~Y.~Hsiang: Cohomology theory of topological
transformation groups, Ergeb. Math. {\bf 85}, Springer, 1975. 

\bibitem{Se} J.-P.~Serre: Valeurs propres des endomorphismes
de Frobenius, S\'eminaire Bourbaki {\bf 446}, 1974. 

\bibitem{S} T.~A.~Springer: Linear algebraic groups. Second edition,
Progress in Math. {\bf 9}, Birkh\"auser 1998.

\end{thebibliography}
\end{document}